\pdfoutput=1
\RequirePackage{ifpdf}
\ifpdf 
\documentclass[pdftex]{sigma}
\else
\documentclass{sigma}
\fi

\usepackage{enumerate}

\numberwithin{equation}{section} 

\newtheorem{prop-defn}{Proposition-Definition}

\DeclareMathOperator{\End}{End}   
\DeclareMathOperator{\Id}{Id}     
\DeclareMathOperator{\id}{Id}     
\DeclareMathOperator{\Tr}{Tr}     
\DeclareMathOperator{\tr}{Tr}     
\DeclareMathOperator{\Adj}{Ad}
\newcommand{\sA}{\mathsf A} 
\newcommand{\sU}{\mathsf U} 
\newcommand{\A}{\mathcal{A}}  
\newcommand{\B}{\mathcal{B}}  
\newcommand{\C}{\mathbb{C}}   
\renewcommand{\H}{\mathcal{H}}  
\newcommand{\J}{\mathcal{J}} 
\newcommand{\N}{\mathbb{N}}   
\newcommand{\R}{\mathbb{R}}   
\newcommand{\Z}{\mathbb{Z}}   
\newcommand{\ox}{\otimes}     
\newcommand{\twobytwo}[4]{\begin{pmatrix}#1 & #2 \\
                           #3 & #4\end{pmatrix}} 
\newcommand{\tst}[3]{{#1}\!\otimes\!{#2}\!\otimes\!{#3}} 

\newcommand{\Ch}{\mathrm{Ch}}
\newcommand{\ket}[1]{\left | #1 \right \rangle}
\newcommand{\ketua}[1]{\ket{ #1 \uparrow}}
\newcommand{\ketda}[1]{\ket{ #1 \downarrow}}

\newcommand{\Wua}{W^\uparrow}
\newcommand{\Wda}{W^\downarrow}
\newcommand{\Pua}{P^\uparrow}
\newcommand{\Pda}{P^\downarrow}
\newcommand{\qn}[1]{\left [ #1 \right ]}
\DeclareMathOperator{\qInd}{q-Ind}
\DeclareMathOperator{\dom}{dom}
\DeclareMathOperator{\coker}{coker}
\DeclareMathOperator{\sign}{sign}

\begin{document}

\allowdisplaybreaks

\renewcommand{\PaperNumber}{051}

\FirstPageHeading

\ShortArticleName{Twisted Cyclic Cohomology and Modular Fredholm Modules}

\ArticleName{Twisted Cyclic Cohomology\\ and Modular Fredholm Modules}

\Author{Adam RENNIE~$^{\dag^1}$, Andrzej SITARZ~$^{\dag^2\dag^3}$ and Makoto YAMASHITA~$^{\dag^4}$}

\AuthorNameForHeading{A.~Rennie, A.~Sitarz and M.~Yamashita}

\Address{$^{\dag^1}$~School of Mathematics and Applied Statistics,  University of Wollongong,\\
\hphantom{$^{\dag^1}$}~Wollongong NSW 2522, Australia}
\EmailDD{\href{mailto:renniea@uow.edu.au}{renniea@uow.edu.au}}

\Address{$^{\dag^2}$~Institute of Mathematics of the
Polish Academy of Sciences, \\
\hphantom{$^{\dag^2}$}~ul.~Sniadeckich 8, Warszawa, 00-950 Poland}

\Address{$^{\dag^3}$~Institute of Physics, Jagiellonian University, ul.~Reymonta 4, 30-059 Krak\'ow, Poland}
\EmailDD{\href{mailto:andrzej.sitarz@uj.edu.pl}{andrzej.sitarz@uj.edu.pl}}

\Address{$^{\dag^4}$~Department of Mathematics, Ochanomizu University, Otsuka 2-1-1, Tokyo, Japan}
\EmailDD{\href{mailto:yamashita.makoto@ocha.ac.jp}{yamashita.makoto@ocha.ac.jp}}

\ArticleDates{Received January 24, 2013, in f\/inal form July 22, 2013; Published online July 30, 2013}

\vspace{-1mm}

\Abstract{Connes and Cuntz showed in [\textit{Comm.\ Math.\ Phys.}\ \textbf{114} (1988), 515--526] that suitable cyclic cocycles
can be represented as Chern characters of f\/initely summable semif\/inite
Fredholm modules. We show an analogous result in twisted cyclic cohomology
using Chern characters of modular Fredholm modules. We present examples
of modular Fredholm modules arising from Podle\'s spheres and from
${\rm SU}_q(2)$.}

\Keywords{twisted cyclic cohomology; spectral triple; modular theory; KMS weight}

\Classification{58J42; 58B32; 46L87}

\vspace{-2.5mm}

\section{Introduction}
\label{sec:intro}

Let $\A$ be an associative algebra over the f\/ield of complex numbers $\C$,
$\A*_\C \A$ the free product, and\footnote{There is a clash
in the standard notations: in this section $q$ is used for $q\A$
and~$q(a)$, while later in the paper $q$~is used as a deformation
parameter. The dif\/ferent usages will be clear from context.} $q\A$ the ideal generated
by $\iota(a)-\bar{\iota}(a)$, $a\in \A$, where $\iota$, $\bar{\iota}$ are the
two canonical inclusions of $\A$ in $\A*_\C \A$. In~\cite{CC},
it was shown  that those cyclic cocycles for $\A$ which arise from
positive traces on~$(qA)^n$ are Chern characters of f\/initely summable
semif\/inite Fredholm modules.

In this note we show that those twisted cyclic cocycles arising from KMS
weights on $(qA)^n$ are Chern characters of f\/initely summable modular
Fredholm modules, a twisted version of the usual notion of Fredholm modules.
While this is not in any
way a practical method of obtaining such representing Fredholm modules, it shows that in
general one must consider the semif\/inite and modular settings.

The examples treated in the last two sections, the Podle\'{s} spheres $S^2_{q, s}$ and the ${\rm SU}_q(2)$ quantum group, show that we can construct non-trivial twisted cyclic cocycles from
naturally arising modular Fredholm modules. Moreover these cocycles encode the
correct classical dimension, in the sense that the Hochschild class of these cocycles
is non-vanishing at the classical dimension. This was f\/irst observed in~\cite{NT05} for the standard Podle\'{s} sphere.  Thus using twisted cohomology avoids the `dimension drop' phenomena, at least in these examples.

We observe that the objects and phenomena studied here seem to have little to do with~\cite{CM08} and related papers such as~\cite{FK10}. The use of twisted commutators in these
papers leads to a need for twisted traces, but ultimately these produce actual (not twisted)
cyclic cocycles.

\section{The algebraic background}
We begin with a short recollection of the  twisted cyclic cohomology of an algebra $\A$.

\begin{definition}
Let $\A$ be an algebra and $\sigma$ be an automorphism of $\A$. We say that
$\phi \colon \A^{\otimes (n+1)} \to \C$ is a $\sigma$-twisted cyclic $n$-cocycle if,
\begin{itemize}\itemsep=0pt
\item $\phi$ is $\sigma$-invariant:
\[
\phi(a_0,a_1,\ldots,a_n) =  \phi(\sigma(a_0),\sigma(a_1),\ldots,\sigma(a_n));
\]
\item $\phi$ is $\sigma$-cyclic:
\[
\phi(a_0,a_1,\ldots,a_n) = (-1)^n \phi(\sigma(a_n),a_0,a_1,\ldots,a_{n-1});
\]
\item $\phi$ is a $\sigma$-twisted Hochschild cocycle
\begin{gather*}
(b^\sigma \phi) (a_0,a_1,\ldots,a_n, a_{n+1}) =  \sum_{k=0}^n (-1)^k \phi(a_0,\ldots, a_k a_{k+1}, \ldots, a_{n+1}) \\
\hphantom{(b^\sigma \phi) (a_0,a_1,\ldots,a_n, a_{n+1}) =}{}
 + (-1)^{n+1}  \phi(\sigma(a_{n+1}) a_0, a_1, \ldots, a_{n}) = 0,
\end{gather*}
\end{itemize}
where $a_0,\ldots,a_{n+1}\in\A$.
\end{definition}

In the examples we consider, one can use
the algebraic tensor product, as we work with polynomial subalgebras. Alternatively,
one could complete $\A$ in a suitable Fr\'echet topology and use the projective tensor product:
see \cite{CGRS2} and \cite[III, Appendix B]{Co} for more information. As all our algebras
will have a natural $C^*$-completion, we def\/ine all $K$-theories in terms of such completions.
It will transpire for our examples that generators of the relevant $K$-theory groups will belong
to the polynomial subalgebras we will work with.

Now, let us present a simple generalisation of a  result of Connes and Cuntz \cite{CC} to the twisted cyclic theory.

Let $\A$ be a unital algebra and $q\A$ be an algebra generated by the elements
$q(a)$, $q(a) b$, and $a q(b)$ for $a, b\in\A$,
subject to the relation $q(\lambda a + b + \mu) = \lambda q(a) + q(b)$ and
\begin{gather}
q(ab) = q(a) b + a q(b) - q(a) q(b), \qquad a, b \in \A.
\label{qA}
\end{gather}
Equivalently, one may identify $q\A$ with the ideal
within the unital free product algebra $\A*_\C\A$ generated by the elements $q(a):=\iota(a)-\bar{\iota}(a)$ for $a\in\A$.
If $\A$ is an involutive algebra, then so is~$q\A$ with the involution def\/ined by $q(a)^* = q(a^*)$ for $a\in\A$.

Setting $\J:=q\A\subset \A*_\C\A$, we can def\/ine $\J^n$ to
be the ideal of $\A*_\C\A$ generated by the products
$a_0 q(a_1) \cdots  q(a_m)$ and $q(a_1) \cdots  q(a_m)$
for $m \geq n$.
If $\sigma$ is an automorphism of $\A$, then we can
extend $\sigma$ to an automorphism of $\J$ and $\J^n$
by setting $\sigma(q(a)):=q(\sigma(a))$.

\begin{proposition}[\protect{see \cite[Proposition 3]{CC}}]
Let $\A$ be a unital algebra, $\sigma$ an automorphism of $\A$, and let $\J$ be the ideal $q\A$ of $\A*_\C\A$ described above.
Suppose that $T$ is a $\sigma$-twisted trace on $\J^n$ for some even integer $n$. That is, $T$ is a
linear functional such that
\begin{gather}
T(xy) = T(\sigma(y)x), \qquad \forall\, x \in J^k,\quad y \in J^l, \quad k+l=n
\label{EqTwist}
\end{gather}
with the convention that $J^0 = \A *_\C \A$.  Then the formula
\[
\tau(a_0,a_1,\dots, a_n) := T(q(a_0) q(a_1) \cdots q(a_n)),\qquad a_0,a_1,\dots,a_n\in\A
\]
defines a $\sigma$-twisted cyclic $n$-cocycle $\tau$ on $\A$.
\end{proposition}

\begin{proof}
Setting $x = 1$ in \eqref{EqTwist}, we obtain that $T$ is $\sigma$-invariant. The $\sigma$-cyclicity follows by setting $x = q(a_0)\cdots q(a_{n-1})$ and $y = q(a_n)$. It remains to verify the $\sigma$-twisted Hochschild cocycle condition.
If $a_0,\dots,a_n$ are elements of $\A$, \eqref{qA} implies
\begin{gather*}
\sum_{k=0}^n (-1)^k T \left( q(a_0) \cdots q(a_k a_{k+1}) \cdots q(a_{n+1}) \right) \\
\qquad{} = T(a_0 q(a_1) \cdots q(a_{n+1})) + T(q(a_0) \cdots q(a_n) a_{n+1}) - T(q(a_0) \cdots q(a_{n+1})).
\end{gather*}
Then, using \eqref{EqTwist}, one sees that this is equal to
\[
T((q(\sigma(a_{n+1})) a_0 + \sigma(a_{n+1}) q(a_0) - q(\sigma(a_{n+1})) q(a_0)) q(a_1) \cdots q(a_n)).
\]
Again by \eqref{qA}, we obtain
\[
\sum_{k=0}^n (-1)^k T \left( q(a_0) \cdots q(a_k a_{k+1}) \cdots q(a_{n+1}) \right) = T \left( q(\sigma(a_{n+1}) a_0)q(a_1) \cdots q(a_{n}) \right),
\]
which is equivalent to the desired equality $b^\sigma \tau = 0$.
\end{proof}

For the analogous statement for odd cocycles, we need to extend the automorphism $\sigma$ to $\J^n$ in a dif\/ferent way, cf.~\cite[Lemma~4]{CC}.
We def\/ine $\tilde{\sigma}$ via the formula
\begin{gather*}
  \tilde{\sigma}\left( a_0 q(a_1) \cdots  q(a_m) \right) =
(-1)^m (\sigma(a_0) - q(\sigma(a_0))) q(\sigma(a_1)) \cdots  q(\sigma(a_m)), \\
 \tilde{\sigma} \left( q(a_1) \cdots  q(a_m) \right)    =
(-1)^m q(\sigma(a_1)) \cdots  q(\sigma(a_m)).
\end{gather*}
Then it is easy to check that $\tilde{\sigma}$ is indeed an
automorphism of $q\A$ and, just as above,  we have

\begin{proposition}
If $T$ is a $\tilde{\sigma}$-twisted trace on~$\J^{n}$, for $n$ an odd integer,
then the formula
\[
\tau(a_0,a_1,\dots, a_{n}) := T(q(a_0) q(a_1) \cdots q(a_{n})),\qquad a_0,a_1,\dots,a_n\in\A,
\]
defines a  $\sigma$-twisted $n$-cyclic cocycle on $\A$.
\end{proposition}

\section{The analytic picture}

In this section we look at a version of \cite[Theorem~15]{CC} in twisted
cyclic cohomology. In brief, \cite{CC} shows that positive traces on certain
ideals in the free product~$\A*_\C\A$ give rise to cyclic cocycles on~$\A$.
These cyclic cocycles can be represented as the Chern characters of semif\/inite
Fredholm modules. By replacing traces with KMS functionals, we arrive at an
analogue of this result in twisted cyclic theory. There are also some analytic
dif\/ferences in our starting assumptions, which we discuss at the end of this section.

We let $\sA$ be a unital $C^*$-algebra and consider the unital
full free product $C^*$-algebra
$\sA*_\C \sA$. We denote by $\iota$, $\bar{\iota}$ the two canonical inclusions
 of $\sA$ in $\sA*_\C \sA$,
and by $q\sA$ the ideal generated by
elements of the form $q(a):=\iota(a)-\bar{\iota}(a)$ for $a \in \sA$.

Similarly, if $\A \subset \sA$ is a dense subalgebra, then
we let $q\A$ be the analogously def\/ined ideal in $\A*_\C\A$.
Introduce the shorthand $J^k:=(q\sA)^k$ and $\J^k:=(q\A)^k$ for $k\in\N$.

Our starting point is a lower semicontinuous and semif\/inite
weight $\phi$ on the $C^*$-algebra~$J^{2p}$ \cite[Chapter~VI]{TT}
which satisf\/ies the KMS$_\beta$ condition for a strongly continuous
one parameter group $\sigma_\bullet\colon \mathbb{R} \rightarrow \mbox{Aut}\big(J^{2p}\big)$.
We will assume that
$\J^{2p}\subset \dom(\phi)$ and that $\J^{2p}$ consists of analytic
vectors for $\sigma_\bullet$, and that
\begin{gather}
\label{EqFaithfulQuot}
\phi(x x^*) = 0 \ \Leftrightarrow \ \phi(x^* x) = 0, \qquad x \in \sA.
\end{gather}

The weight $\phi$ gives, via the GNS construction, a Hilbert space $\H_\phi$ with
a nondegenerate representation $\pi_\phi\colon J^{2p}\to \B(\H_\phi)$, and a linear map
$\Lambda\colon \dom^{1/2}(\phi)\subset J^{2p}\to \H_\phi$. The condition~\eqref{EqFaithfulQuot} implies that $\{ x \in \sA \mid \phi(x x^*) = 0 \}$ is the kernel of this representation, and that $\pi_\phi\big(J^{2p}\big)$ admits a faithful weight which induces $\phi$.  Since $\sigma_\bullet$ leaves $\phi$ invariant, it descends to $\pi_\phi(J^{2p})$.

There is a canonical faithful normal semif\/inite extension $\Phi$ of
$\phi$ to $\big(\pi_\phi\big(J^{2p}\big)\big)''$ satisfying $\phi=\Phi\circ \pi_\phi$
and $\sigma^\Phi_t\circ \pi_\phi=\pi_\phi\circ \sigma_{-\beta t}$. See \cite[Proposition~1.5, Chapter~VIII]{TT} for
a proof.

The KMS property implies that for $a,\,b\in J^{2p}$ we have
\[
\phi(ab)=\phi(\sigma(b)a),
\]
where we def\/ine the (non-$*$) automorphism $\sigma$ to be the value of the extension of the one-parameter group
$\sigma_\bullet$ to the complex value $t=i\beta$.

We observe that the representation of $J^{2p}$ on $\H_\phi$ extends naturally to
a representation of $\sA*_\C \sA$ on $\H_\phi$, denoted $\lambda$,
such that $\lambda(\sA*_\C \sA)\subset \big(\pi_\phi\big(J^{2p}\big)\big)''$.  This is the usual
extension, def\/ined on the dense subspace $\pi_\phi\big(J^{2p}\big)\H_\phi$
by $\lambda(\alpha)(j\xi):=(\alpha j)\xi$ for $\alpha\in \sA*_\C \sA$, $j\in J^{2p}$
and $\xi\in \H_\phi$. If~$T$ is in the commutant of
$\pi_\phi(J^{2p})$ then
\[
T(\lambda(\alpha)(j\xi))=T((\alpha j)\xi)=(\alpha j)(T\xi)=\lambda(\alpha)(j(T\xi))=\lambda(\alpha)(T(j\xi)),
\]
showing that $\lambda(\alpha)$ is indeed in $\pi_\phi(J^{2p})''$ for all
$\alpha\in \sA*_\C \sA$.

By \cite[Theorem~2.6, Chapter~VII]{TT}, the (image under $\Lambda$ of)
$\dom^{1/2} \Phi\cap\big(\dom^{1/2} \Phi\big)^*$ is a full left Hilbert
algebra, which we denote by $\sU$. Moreover, the left von Neumann algebra of
$\sU$ is precisely $\big(\pi_\phi\big(J^{2p}\big)\big)''$. We record the following Lemma,
whose proof follows immediately from the def\/initions.

\begin{lemma}
Let $N$ be the left von Neumann algebra of the left Hilbert algebra~$U$
and~$\Phi$ the corresponding faithful normal semifinite weight. Then
for all $\alpha\in J^p \cap \dom^{1/2}(\phi)$ we have
$\lambda(\alpha)\in \dom^{1/2}(\Phi)$, and $\Phi(\lambda(\alpha)^*\lambda(\alpha))=\phi(\alpha^*\alpha)$.
\end{lemma}

\begin{definition}
Let ${\mathcal N}$ be a von Neumann algebra acting on a Hilbert space $\H$, and $\Phi$ be a~faithful normal
semif\/inite weight on $\mathcal{N}$.
Then we say that $(\A,\H,F)$ is an $n$-summable {\it unital modular Fredholm module}
with respect to $({\mathcal N},\Phi)$ if
\begin{enumerate}[i)]\itemsep=0pt
\item[o)] $\A$ is a separable unital $*$-subalgebra of ${\mathcal N}$;
\item   $\A$ is globally invariant under the group $\sigma^\Phi$, and consists of
             analytic vectors for it;
\item   $F$ is a self-adjoint operator in the f\/ixed point algebra
            ${\mathcal M}:={\mathcal N}^{\sigma^\Phi}$ with $F^2=1_{\mathcal N}$;
\item  $[F,a]^{n}\in \dom(\Phi)$ for all $a\in \A$.
\end{enumerate}
If there exists a self adjoint element $\gamma$ of ${\mathcal M}$ satisfying $\gamma^2=1$, $\gamma a=a\gamma$ for all $a\in\A$ and $\gamma F + F \gamma=0$, the quadruple $(\A, \H, F, \gamma)$ is said to be an even module.
In contrast, a Fredholm module without the grading $\gamma$ is said to be odd.
\end{definition}

The Chern character of a modular Fredholm module is the class of
the $\sigma:=\sigma^\Phi_i$ twisted cyclic $n$-cocycle def\/ined by
the formula
\[
\hbox{Ch}_{n}(a_0,a_1,\dots,a_{n})=\lambda_n\frac{1}{2}\Phi(\gamma F[F,a_0][F,a_1]\cdots[F,a_n]),\qquad a_0,a_1,\dots,a_n\in\A.
\]
Here we set $\gamma=1$ if the module is odd. The constants $\lambda_n$ are given by
\[
\lambda_n= \begin{cases} \displaystyle (-1)^{n(n-1)/2}\Gamma\left(\frac{n}{2}+1\right), & n\ \mbox{even},\vspace{1mm}\\
\displaystyle \sqrt{2 i}(-1)^{n(n-1)/2}\Gamma\left(\frac{n}{2}+1\right), & n\ \mbox{odd}.\end{cases}
\]

\begin{theorem}
\label{thm:main-twist}
Suppose that $\A$ is a $*$-subalgebra of a $C^*$-algebra $\sA$,
and $\phi$ is a weight on $J^{2p}$ which is lower semicontinuous, semifinite,
and satisfies~\eqref{EqFaithfulQuot}.  We further assume
that it satisfies the KMS$_\beta$ condition for a one parameter
group $\sigma$ such that $\J^{2p}:=(q\A)^{2p}$ consists of analytic
vectors in the domain of $\phi$.
Then there exists a $2p$-summable modular Fredholm module for
$\A$. The modular
Fredholm module has Chern character
\[
{\rm Ch}_{2p}(a_0,a_1,\dots,a_{2p})=\lambda_{2p}(-1)^p\phi(q(a_0)\,q(a_1)\cdots q(a_{2p})).
\]
\end{theorem}

\begin{proof}
The universal property of $\sA*_\C \sA$ gives two $*$-homomorphisms
$\pi_\phi$ and $\bar{\pi}_\phi$ from $\sA$ to $\B\big(L^2\big(J^{p+1},\phi\big)\big)$, whose images lie
in $N=\pi_\phi\big(J^{2p}\big)''$.
The modular Fredholm module is given  by the data:
\begin{itemize}\itemsep=0pt
\item the Hilbert space $\H:=L^2\big(J^{p+1},\phi\big)\oplus L^2\big(J^{p+1},\phi\big)$;
\item the representation $\pi_2\colon \A\to \B(\H)$, $\pi_2(a)=\pi_\phi(a)\oplus\bar{\pi}_\phi(a)$;
\item the operator $F=\begin{pmatrix} 0 & 1\\ 1 & 0\end{pmatrix}$;
\item the von Neumann algebra $M_2(N)$;
\item the weight $\Phi\circ \mbox{Tr}_2$.
\end{itemize}
Observe that
\[
[F,\pi_2(a)]=\begin{pmatrix} 0 & \bar{\pi}_\phi(a)-\pi_\phi(a)\\ \pi_\phi(a)-\bar{\pi}_\phi(a) & 0\end{pmatrix}=\begin{pmatrix} 0 & -\pi_\phi(q(a))\\ \pi_\phi(q(a)) & 0\end{pmatrix}.
\]
Since $J^{2p+1}\subset J^{2p}$, the Chern character is well-def\/ined. The computation
of the Chern character is straightforward.
\end{proof}

The odd case is similar using the
$\tilde{\sigma}$-automorphism of~$\J^n$.

\begin{remark}
When $\sigma$ is trivial, the above construction reduces to a particular case of the one in \cite[Section~V]{CC}.  Since we started from a positive (unbounded) functional on the enveloping $C^*$-algebra rather than a one on the algebraic object~$\J^{2 p}$, the analytic argument leading to the existence of the underlying left Hilbert algebra is greatly simplif\/ied.  If we were to assume only
a `positive twisted trace' on~$\J^{2p}$ in some sense, we would not necessarily
have enough control to prove this pre-closedness.
\end{remark}

\section[The modular index and pairing with $K$-theory]{The modular index and pairing with $\boldsymbol{K}$-theory}

We recall here the construction of the modular index and its computation
through the pairing between the equivariant $K$-theory and twisted cyclic
cohomology. This section adapts \cite{NT04} to the notation and notions
used here.

Let ${\mathcal N}$ be a von Neumann algebra endowed with a faithful normal semif\/inite weight $\Phi$, and $(\A, \H, F, \gamma)$ be a $2n$-summable even modular Fredholm module
with respect to $({\mathcal N},\Phi)$, as def\/ined in the previous section.

Furthermore, let us assume that there exists a densely def\/ined operator $\Xi $
in $\H$ which implements the modular automorphism.  Thus, $\Xi$ is an unbounded self adjoint operator satisfying
\[
 [F,\Xi ]=0, \qquad [\Xi ,\gamma]=0,\qquad \Xi ^{-1} a \Xi  = \sigma(a),
\qquad \forall \, a \in \A,
\]
where we identify $a$ with the operator in the representation $\pi(\A)$.

With this set up we can make the following def\/inition, and we assume in what follows that
$({\mathcal N},\Phi)=\big(\B(\H),\tr\big(\Xi^{1/2}\cdot\Xi^{1/2}\big)\big)$, as this is the context we shall
be working with in the examples. Extending the def\/inition to the more general situation
is straightforward using the theory of Breuer--Fredholm operators as in~\cite{CPRS3}.

\begin{definition}\label{defn:even-modular-index}
Let $F$ be a Fredholm operator, which commutes with $\Xi$. We def\/ine the \textit{modular
index} of $F$ to be
\[
\qInd(F) = \tr (\Xi |_{\ker F}) - \tr (\Xi |_{\coker F}).
\]
\end{definition}

This def\/inition is well def\/ined, since both kernel and cokernel are f\/inite-dimensional
and at the same time invariant subspaces of $\Xi$, so in fact both traces are f\/inite
expressions.  We omit the proof of the next standard construction.

\begin{prop-defn}
Suppose that $(\A, \H, F, \gamma)$ is an even modular Fredholm module, and let $p \in \sA$ be a projection which is fixed by the
modular automorphism group~$\sigma_\bullet$.  Replacing~$\H$ by~$p \H$, $\mathcal{N}$ by $\mathcal{N}_p = p \mathcal{N} p$, $\Phi$ by $\Phi_p = \Phi|_{\mathcal{N}_p}$, $F$ by $p F_+ p$ in Definition~{\rm \ref{defn:even-modular-index}}, we obtain a
Fredholm operator for $(\mathcal{N}_p, \Phi_p)$.  We define
$\qInd^F(p)$ to be its modular index $\qInd(p F_+ p)$.
\end{prop-defn}

More generally, we can extend the above index pairing on the classes in the equivariant $K_0$-group as follows.  An element of the
equivariant $K_0$-group $K_0^\R(\A)$ is given by a formal
dif\/ference of invariant projections in the $\R$-algebras of the
form $\A \otimes \End(X)$, where $U\colon \R\to\End(X)$ is an arbitrary f\/inite-dimensional representation of $\R$~\cite[Theorem~3.1]{NT04}.  Assume that
$p \in \A \otimes \End(X)$ is such a projection.
Then we extend the modular Fredholm module to
$(\A \otimes \End(X), \H \otimes X,
F \otimes \id, \gamma \otimes \id)$ with respect to
$(\mathcal{N} \otimes \End(X), \Phi \otimes G_X)$, where $G_X(T)=\tr(U_{-i}T)$ for
$T\in\End(X)$.
The above consideration gives the number $\qInd^F(p)$,
which only depends on the $K_0^\R$-class of $p$, \cite[Lemma 3.15]{RS}.
This way, we obtain the map
\[
\qInd^F\colon K_0^\R(\A) \to \R.
\]

With this def\/inition the following two propositions follow as in \cite{NT04}.

\begin{proposition}\label{Prop:even-mod-fred-mod-ind-pairing}
Let $\H = \H_+ \oplus \H_-$, where $\gamma_{\H_\pm} = \pm \id$ and let us
denote $F_+$ the restriction of $F$ to $\H_+$. Let $U\colon \R\to \End(X)$ be a unitary representation of $\R$ on a finite-dimensional Hilbert space $X$.  Let $D$ be the generator of $U$, and put $\Delta_X = e^{-D}$.  For any projection $p \in \A \otimes \End(X)$
invariant under the action $(\sigma_t \otimes \Adj_{\Delta_X^{i t}})_{t \in \R}$,
the modular index of $p (F_+ \otimes \Id_X) p$ is given by
\begin{gather}
\qInd^F(p) = (-1)^n \tr_{\H} \otimes \tr_X ((\Xi \otimes \Delta_X) (\gamma \otimes \Id_X)
p [F\otimes \Id_X, p]^{2n}).
\label{qind-even}
\end{gather}
\end{proposition}

\begin{proof}[Outline of the proof]
The proof follows standard lines, see \cite[pp.~296--297]{Co} and \cite[pp.~370--371]{NT04}.
The broad plan is to observe that
\[
p-p(F\ox\Id_X)p(F\ox\Id_X)p=-p[(F\ox\Id_X),p][(F\ox\Id_X),p]p.
\]
Then we observe that
$\tau:=\tr_{\H} \otimes \tr_X ((\Xi \otimes \Delta_X)\cdot)$ is a trace on the f\/ixed point algebra
$(\B(\H)\ox\End(X))^{\sigma\ox {\rm Ad}\Delta_X}$.  Since
$(p-p(F\ox\Id_X)p(F\ox\Id_X)p)^n$ is trace class with respect to $\tau$, we
f\/ind that the $\tau$-index of $p(F_+\ox\Id_X)p$ is given by
$\tau(\gamma(p-p(F\ox\Id_X)p(F\ox\Id_X)p)^n)$ (see \cite[Proposition 4.2]{GVF}). By
def\/inition of $\tau$, this is just
$\qInd^F(p)$.
\end{proof}

There is the notion of modular index pairing for odd modular Fredholm modules and invariant unitaries.  Let $(\A, \H, F)$ be an odd modular Fredholm module with respect to $(\mathcal{N}, \Phi)$, and set $E=\frac{1}{2}(1+F)\otimes \Id_X$.

\begin{prop-defn}
Let $X$ be a finite-dimensional unitary representation of $\R$, and suppose that $v$ is an  unitary element in $\A\ox\End(X)$ which is invariant under $(\sigma_t \otimes \Adj_{\Delta_X^{i t}})_{t \in \R}$.
Then, $E v E v^*$, as an operator from $v E (\H\ox X)$ to $E(\H\ox X)$, becomes a modular Fredholm ope\-ra\-tor.  We let $\qInd^F(v)$ denote its modular index.  This number only depends on the equivariant $K_1$-class of $v$.
\end{prop-defn}

\begin{proposition}
Under the above setting, the modular index of $E v E v^*$ is given by
\begin{gather}
\qInd^F(v) = \frac{(-1)^n}{2^{2n}} \tr_{\H}\otimes\tr_X \big(\Xi\otimes\Delta_X
\left(  [F,v] [F,v^{*}] \right)^n\big).
\label{qind-odd}
\end{gather}
\end{proposition}

\begin{proof}[Outline of the proof]
As in the even case, we observe that $(E-vEv^*)^2=-\frac{1}{4}E[F,v][F,v^*]$, and so by
\cite[Theorem~3.1]{CP2} and the def\/inition of $\tau$ above,
the $\tau$-index of $EvEv^*$ is given by
$\tau((E-vEv^*)^{2n}-(E-v^*Ev)^{2n})$, and after standard algebraic manipulations  (see~\cite[pp.~51--52]{CGRS2}), we f\/ind the result of the proposition.
\end{proof}

\section{The Podle\'s spheres}

\subsection{The algebra}

Given parameters $0 \leq q < 1$ and $0 \leq s \leq 1$, the Podle\'s quantum sphere $\A\big(S^2_{q,s}\big)$ is def\/ined as the universal $*$-algebra with generators $A=A^*$, $B$, and $B^*$ subject to the relations
\begin{gather*}
B^*B + (A-1)\big(A+s^2\big)  = 0, \qquad
BB^* + \big(q^2A-1\big)\big(q^2A+s^2\big)  = 0, \qquad
AB  = q^{-2} BA.
\end{gather*}

When $0 < s$, the algebra $\A\big(S^2_{q,s}\big)$ has two inequivalent irreducible representations $\pi_+$ and $\pi_-$ on $\ell^2(\N)$. In terms of is the standard orthonormal basis $\{e_k\}_{k\in\N}$, they are given by the formulae
\begin{alignat}{3}
& \pi_+(B) e_k= \sqrt{1-q^{2k}}\sqrt{s^2+q^{2k}} e_{k-1}, \qquad &&
\pi_+(A) e_k  =  q^{2k} e_k, &\nonumber\\
& \pi_-(B) e_k= s \sqrt{1-q^{2k}}\sqrt{1+s^2 q^{2k}} e_{k-1}, \qquad & &
\pi_-(A) e_k  = -s^2 q^{2k} e_k. &\label{eq:podles-sph-action}
\end{alignat}

The algebra $\A\big(S^2_{q,s}\big)$ can be completed to a $C^*$-algebra $C\big(S^2_{q,s}\big)$ by means of the
operator norm induced by the representation $\pi_+ \oplus \pi_-$.
The modular group $(\sigma_t)_{t \in \R}$ is periodic, and it extends by continuity to an action of $\mathrm{U}(1)$ on $C\big(S^2_{q,s}\big)$.  These $\mathrm{U}(1)$-$C^*$-algebras are
known to be isomorphic to the f\/ibre product of the two copies of the Toeplitz algebra $\mathcal{T}$ with
respect to the symbol map $\mathcal{T} \rightarrow C(S^1)$~\cite{S}.  Here,
we consider the gauge action of $\mathrm{U}(1)$ on each copy of $\mathcal{T}$ and the translation action on $C(S^1)$.

\begin{definition}
We construct an even Fredholm module $\big(\A\big(S^2_{q,s}\big), F, \H\big)$ by taking
$\H = \ell^2(\N) \oplus \ell^2(\N)$, endowed with the representation of $\A\big(S^2_{q,s}\big)$ def\/ined by the formula
\[
\pi(a) = \left( \begin{matrix} \pi_+(a) & 0 \\ 0 & \pi_-(a) \end{matrix} \right), \qquad a \in \A\big(S^2_{q,s}\big),
\]
along with the grading operator and Fredholm operator $F$ given by
\begin{gather*}
\gamma  = \left( \begin{matrix} 1 & 0 \\ 0 & -1 \end{matrix} \right), \qquad
F  = \left( \begin{matrix} 0 & 1 \\ 1 & 0 \end{matrix} \right).
\end{gather*}
In addition, we let $K$ be the diagonal modular operator on $\H$ def\/ined by
\[
K e_{k,\pm} = q^{-2k} e_{k,\pm},\qquad k \in \N.
\]
Here, $e_{k, +}$ is a basis in the f\/irst direct summand (supporting $\pi_+$) $\ell^2(\N)$ of~$\H$, and $e_{k, -}$ is in the second direct summand (supporting $\pi_-$).
\end{definition}

\begin{lemma}\label{lemma5.1}
For $0<s\leq 1$ and $0<q<1$, the Fredholm module $\big(\A\big(S^2_{q,s}\big), F, \H\big)$ can be regarded as a $2$-summable modular Fredholm module with the following data.  The von Neumann algebra~$\mathcal{N}$ is $\B(\H)$, with the weight $\Phi$ defined by
\[
\Phi(S) = {\rm Tr}\big(K^{1/2} S K^{1/2}\big),\qquad 0 \leq S \in \B(\H).
\]
The modular automorphism $\sigma(T) = K^{-1} T K$ leaves $\A\big(S^2_{q,s}\big)$ invariant, and the restriction can be expressed as
\begin{gather*}
\sigma(A)  = A,\qquad
\sigma(B)  = q^{-2} B,\qquad
\sigma(B^*)  = q^2 B^*.
\end{gather*}
\end{lemma}

\begin{proof}
By \cite[Theorem 2.11]{TT}, $\Phi$ is a faithful normal semif\/inite weight on $\B(\H)$, with modular
group given by $T\mapsto K^{it}TK^{-it}$.

As a next step, we show that for any $a \in \A\big(S^2_{q,s}\big)$ the operator
$K [F,\pi(a)]$ is bounded and $[F,\pi(a)]$ is of trace class. Applying the def\/initions
of the representation $\pi$ and the operator $F$ yields
\begin{gather*}
  K [F,\pi(A)] e_{k,\pm} = \pm \big(1+s^2\big) e_{k,\mp}, \\
  K [F,\pi(B)] e_{k,\pm} = \pm q^{-2k} \sqrt{1-q^{2k}}
   \left( \sqrt{s^2 + q^{2k}} - s \sqrt{1+ s^2 q^{2k}} \right) e_{k-1,\mp}.
\end{gather*}
To make an estimate for the last expression we observe that
\begin{gather*}
\sqrt{s^2 + q^{2k}} - s \sqrt{1+ s^2 q^{2k}} =
 \frac{ \big(s^2 + q^{2k}\big) - s^2 \big(1+ s^2 q^{2k}\big) }{ \sqrt{s^2 + q^{2k}} + s \sqrt{1+ s^2 q^{2k}} }
 = \frac{ \big(1 - s^4\big) q^{2k} }{ \sqrt{ s^2 + q^{2k}} + s \sqrt{1+ s^2 q^{2k}} }.
\end{gather*}
Since the denominator is greater than or equal to $2s$, we f\/ind
\[
\left| \sqrt{s^2+q^{2k}} - s \sqrt{1+ s^2 q^{2k}} \right| =
\frac{ \big(1 - s^4\big) q^{2k} }{ \sqrt{s^2 + q^{2k}} + s \sqrt{1+ s^2 q^{2k}} } \leq
\frac{1 - s^4}{2s} q^{2k}.
\]

Now, since $[F,\pi(a)] = K^{-1} (K [F,\pi(a)])$ and $K^{-1}$ is a trace
class operator, it follows directly that $[F,\pi(a)]$ is of trace class for any
$a \in \A(S^2_{q, s})$.

Therefore, in the end we obtain that for any $a_0,a_1 \in  \A\big(S^2_{q,s}\big)$,
the operator $ K [F,a_0][F, a_1]$ is of trace class, and so $[F,a_0][F,a_1]$ is
in the domain of $\Phi$.
\end{proof}

\begin{corollary}
The $3$-linear functional $\phi$ defined by the formula
\[
\phi(a_0,a_1,a_2) = \Phi(\gamma F [F,a_0] [F, a_1] [F, a_2]),
\qquad a_0,a_1,a_2\in \A\big(S^2_{q,s}\big),
\]
determines a $\sigma$-twisted cyclic cocycle over $\A\big(S^2_{q,s}\big)$.
\end{corollary}

To see that the cyclic cocycle we obtained is non-trivial,
we explicitly compute  its pairing with the twisted cyclic cycle
$\omega_2$, found by Hadf\/ield \cite{H}. In our notation the twisted cyclic 2-cycle $\omega_2$
is given by
\begin{gather*}
\omega_2 =  2 \left( \tst{A}{B}{B^*} -\tst{A}{B^*}{B} + 2 \tst{B}{B^*}{A}
           -2 q^{-2} \tst{B}{A}{B^*} \right. \\
\left.\hphantom{\omega_2 =}{}       + \big(q^4-1\big) \tst{A}{A}{A} \right) +\big(1-q^{-2}\big) s^2\big(1-s^2\big) \tst{1}{1}{1} \\
\hphantom{\omega_2 =}{}
       + \big(1-s^2\big) \left(\tst{1}{B^*}{B}-q^{-2} \tst{1}{B}{B^*}+\big(1-q^2\big) \tst{1}{A}{A } \right).
\end{gather*}

The pairing of $\omega_2$ with $\phi$
(we skip the straightforward computations) gives
\[
\left(\phi, \omega_2 \right) =   \big(1+s^2\big)^3.
\]

Since $\omega_2$ comes from ${\rm HH}_2^\sigma\big(\A\big(S^2_{q,s}\big)\big)$, this also shows that the
Hochschild class of $\phi$ is nontrivial.

\subsection{The index pairing, local picture}

Consider the projection $P \in M_2\big(\A\big(S^2_{q,s}\big)\big)$ def\/ined by
\[
P = \frac{1}{1+s^2} \left( \begin{matrix} 1 - q^2 A &  B \\
 B^* & A + s^2 \end{matrix} \right).
\]
This projection becomes $\R$-invariant with respect to the representation $(\sigma_t \otimes \Delta^{i t})_{t \in \mathbb{R}}$ of $\mathbb{R}$ on $\A\big(S^2_{q,s}\big) \otimes \mathbb{C}^2$, where $\Delta \in M_2(\mathbb{C})$ is given by
\begin{gather}
\Delta = \left( \begin{matrix} q^{-1} & 0 \\
                                             0 & q \end{matrix} \right).
\label{eq:rho-K}
\end{gather}

The classes of $1$ and $P$ generate the $\mathrm{U}(1)$-equivariant $K_0$ group of $C\big(S^2_{q,s}\big)$~\cite{Wagner}.  We explicitly compute  the index pairing of $P$ with the twisted cocycle
given by the Chern character of the modular Fredholm module constructed  above.

\begin{proposition}\label{Prop:cher-F-and-P-pairing}
The pairing of the Chern character of $F$ and $[P]$ is equal to $q$.
\end{proposition}

\begin{proof}
Expanding the relevant def\/initions, one has
\begin{gather*}
 {\rm Ch}_2(P,P,P) =  -\frac{1}{2} \frac{2q}{(1+s^2)^2}\\
{}\times \sum_{k=0}^\infty
\left( 2 s^2 \big(q^4-1\big) q^{4k} +  2 s (a_{k+1} - a_{k})  + \big(q^2-1\big)\big(1-s^2\big)^2 q^{2k} + 2 s (b_k - b_{k+1}) \right),
\end{gather*}
where
\begin{gather*}
 a_k  = \sqrt{s^2 + q^{2k}} \sqrt{1 +s^2 q^{2k}} q^{2k},\qquad
 b_k  = \sqrt{s^2+ q^{2k}} \sqrt{1+s^2 q^{2k}}.
\end{gather*}

We compute the sum explicitly. First, observe that
\[
 \sum_{k=0}^\infty (a_{k+1} - a_{k}) = - \big(1+s^2\big),
\qquad   \sum_{k=0}^\infty (b_{k} - b_{k+1}) = \big(1+s^2\big) - s,
\]
which then allows the rest of the sum to be computed to yield
\begin{gather*}
{\rm Ch}_2(P,P,P)  = -\frac{1}{2} \frac{2q}{\big(1+s^2\big)^2}
    \left( - 2s^2  - 2s \big(1+s^2\big)  -\big(1-s^2\big)^2 + 2s \big(1+s^2-s\big) \right) =  q.
\end{gather*}
This proves the assertion.
\end{proof}

Note that the index pairing is independent of the parameter $s$, since the $K$-theoretic data is invariant under continuous deformation as we shall see in detail in the next section.

\subsection{The index pairing, global picture}

In this section we give an alternative global picture of the index pairing.  For this purpose it is convenient to use a dif\/ferent set of generators of the equivariant $K$-theory group.

Since $A\in \A\big(S^2_{q,s}\big)$ is $\sigma$-invariant, the spectral projections of the selfadjoint
operator $A$ give elements of the $\sigma$-equivariant $K_0$ group of $C\big(S^2_{q,s}\big)$.
For each $k \in \N$, let us denote the projection onto the span of $e_{k,+}$ (resp.\ $e_{k,-}$) by $p^{(+)}_k$ (resp.\ by $p^{(-)}_k$).

By~\eqref{eq:podles-sph-action}, the spectral projections of $A$
for the positive (resp.\ negative) eigenvalues are given by the
$p^{(+)}_k\oplus 0$ (resp.\ $0\oplus p^{(-)}_k$) for $k \in \N$.  In what follows we abbreviate these projections as~$p^{(+)}_k$ and~$p^{(-)}_k$.

From the description of $K$ in Lemma~\ref{lemma5.1}, Proposition~\ref{Prop:even-mod-fred-mod-ind-pairing} implies
\[
{\rm Ch}_2\big(p^{(\pm)}_k,p^{(\pm)}_k,p^{(\pm)}_k\big) = \pm q^{-2 k}, \qquad k \in \N.
\]
Let us relate this computation to the calculation of Proposition~\ref{Prop:cher-F-and-P-pairing}.  One may think of $P$ as a~family of projections on $\mathcal{H} \otimes \C^2$ parametrized by $0 \le q < 1$ and $0  < s \le 1$.  Moreover, the $C^*$-algebras $C(S^2_{q, s})$ can be identif\/ied with each other because they have the same image under the representation~$\pi$.
From its presentation, $P$~is operator norm continuous in the parameters $q \in [0, 1)$ and $s \in (0, 1]$.
Hence the class of $P$ in $K_0^{\mathrm{U}(1)}C\big(S^2_{q,s}\big)$ is independent of~$q$ and~$s$.

Now, the projection $(\pi_+ \oplus \pi_-)(P)$ at $q = 0$ and $s = 1$ can be written as
\[
\frac{1}{2} \twobytwo{1}{S^*}{S}{1 + p^{(+)}_0} \oplus \frac{1}{2} \twobytwo{1}{S^*}{S}{1 - p^{(-)}_0},
\]
where $S$ is the isometry $e_k \mapsto e_{k+1}$ on $\ell^2(\N)$.  This projection and
\[
\twobytwo{1}{0}{0}{p^{(+)}_0} \oplus \twobytwo{1}{0}{0}{0}
\]
are connected by a continuous path of $\mathrm{U}(1)$-invariant projections $f^{(+)}_t \oplus f^{(-)}_t$ ($t \in [0, 1]$) def\/ined by
\begin{gather*}
f^{(+)}_t  = \frac{1}{2} \twobytwo{\sqrt{1 - t^2} + 1}{t S^*}{t S}{(1 - \sqrt{1 - t^2})(1 - p^{(+)}_0) + 2 p^{(+)}_0}
\end{gather*}
and
\begin{gather*}
f^{(-)}_t  = \frac{1}{2} \twobytwo{\sqrt{1 - t^2} + 1}{t S^*}{t S}{(1 - \sqrt{1 - t^2})(1 - p^{(-)}_0)}.
\end{gather*}
Using the representation of $K$ in equation \eqref{eq:rho-K}, we obtain
\[
{\rm Ch}_2(P,P,P)  = {\rm Ch}_2(1,1,1) q^{-1} +
{\rm Ch}_2\big(p^{(+)}_0, p^{(+)}_0, p^{(+)}_0\big) q = q,
\]
which gives a `global' picture of the index pairing.

\section[The modular Fredholm modules over $\A({\rm SU}_q(2))$]{The modular Fredholm modules over $\boldsymbol{\A({\rm SU}_q(2))}$}

As an example of odd-dimensional case, let us present the quantum group ${\rm SU}_q(2)$. In this section the parameter $q$ takes value in $(0, 1)$.  The $*$-algebra $\A({\rm SU}_q(2))$ is universally generated by $a$ and $b$
satisfying the relations
\begin{gather*}
  b a = q ab, \qquad bb^* =b^*b, \qquad b^* a = q a b^*,\qquad
  a a^* + b b^* =1, \qquad a^*a + q^2 b b^* =1.
\end{gather*}

In this section we shall demonstrate that the fundamental Fredholm
module presented f\/irst in~\cite{MNW} and the Fredholm module arising
from the spectral triple constructed in~\cite{DLSSV} both give rise to non-trivial twisted
cyclic cocycles. We explicitly compute the pairing of these cocycles with an element from
the equivariant $K_1$ group, and show that the two pairings are both non-zero.

\subsection{The basic Fredholm module}
\label{sec:the-basic-fred}

We brief\/ly review the construction of the module Fredholm module.  The Hilbert space is $\H = \ell^2(\N) \otimes \ell^2(\Z)$, with an representation $\pi_0$ of the $\A({\rm SU}_q(2))$ def\/ined by
\begin{gather*}
\pi_0(a) e_{k,l}  = \sqrt{1-q^{2k+2}} e_{k+1,l}, \qquad
\pi_0(b) e_{k,l}  =  q^{k} e_{k,l+1},
\end{gather*}
in terms of the standard basis for $k \geq0$, and  $l \in \Z$.
The Fredholm operator $F$ is chosen to be $ F e_{k,l} = \sign(l) e_{k,l}$, where we put $\sign(0) = 1$.

\begin{lemma}
The triple $(F,\pi_0,\H)$ is a $3$-summable modular Fredholm module
with respect to the von Neumann algebra $\B(\H)$ and weight $\Phi$ defined as follows.
Define the modular operator by
\[
K e_{k,l} = q^{-2k} e_{k,l}.
\]
Then the weight $\Phi$ is given by $\Phi(T):= \tr\big(K^{1/2}TK^{1/2}\big)$,  for $T\geq 0$.
\end{lemma}

\begin{proof}
From the way $\pi_0$ is def\/ined, one obtains the commutation relation
\begin{gather}\label{eq:MNW-commutation-with-F}
[F, \pi_0(a)]  = 0, \qquad [F, \pi_0(b)] e_{k,l}  = 2 q^k \delta_{l,-1} e_{k,l+1}.
\end{gather}
It follows that for any $x \in \A({\rm SU}_q(2))$, the matrix coef\/f\/icient of $[F, \pi_0(x)]$ decays by the order of $q^{k}$ with respect to the index $k \in \N$.  Therefore, for any elements $x$, $y$ in $\A({\rm SU}_q(2))$, the operator
$ K [F, \pi_0(x)] [F, \pi_0(y)]$ is bounded, and for any three elements $x$, $y$, and $z$, the operator
$ K [F, \pi_0(x)] [F, \pi_0(y)] [F, \pi_0(z)]$ is of trace class.  This means that $[F, \pi_0(\A({\rm SU}_q(2)))]^3$ is in the domain of $\Phi$.
\end{proof}

The modular automorphism $\sigma(T) = K^{-1} T K$ is given on the generators by
\begin{gather*}
\sigma(b)  = b,\qquad  \sigma(b^*)  = b^*, \qquad \sigma(a)  = q^{-2} a, \qquad \sigma(a^*)  = q^2 a^*.
\end{gather*}

Since the Fredholm module is odd, we can use it to construct a twisted
$3$-cyclic cocycle.  Let us investigate the pairing of this cocycle with the equivariant $K_1$-group.  Consider the unitary~$V$ in~$\A({\rm SU}_q(2)) \otimes M_2(\C)$ given by
\begin{gather}
V =  \left( \begin{matrix} -q b^* & a \\
a^* & b \end{matrix} \right).
\label{eqK1}
\end{gather}
This gives a generator of the $U(1)$-equivariant $K_1$-group of $C({\rm SU}_q(2))$ with respect to~$\sigma_t$.
If we extend the action of the modular operator to
$\H \otimes \C^2$ by the generator
\[
\tilde{K} = K \otimes
\left( \begin{matrix} q^{-1} & 0 \\
       0 &  q \end{matrix} \right),
\]
then we see that $V$ is invariant under $\sigma = \mathrm{Ad}_{\tilde{K}}$.  We can now compute the modular index pairing using~(\ref{qind-odd}).

\begin{proposition}
The equivariant pairing between the class of ${\rm Ch}_3$
and the class of $V$ is nontrivial and is equal to $q$.
\end{proposition}

\begin{proof}
We compute the pairing explicitly:
\begin{gather*}
\left<[\hbox{Ch}_3],  [V] \right>   =
\frac{(-1)^2}{2^4} \Tr
\big( \tilde{K} F [F, V] [F, V^*] [F,V] [F, V^*] \big) \\
\hphantom{\left<[\hbox{Ch}_3],  [V] \right>}{}
 = \frac{1}{16} \left( \sum_{k=0}^\infty \sum_{l \in \Z} 16 q^{-2k} \sign(l)
    \big( q \delta_{l,0} q^{4k} + q^3 \delta_{l,-1} q^{4k} \big) \right) \\
\hphantom{\left<[\hbox{Ch}_3],  [V] \right>}{}
 =  \sum_{k=0}^\infty q \big(1-q^2\big) q^{2k} = q. \tag*{\qed}
\end{gather*}
  \renewcommand{\qed}{}
\end{proof}

It follows from results of Hadf\/ield and Kr\"ahmer (see \cite[Lemma~4.6]{HK}) that the map
\[
I\colon  \ {\rm HH}^\sigma_3(\A({\rm SU}_q(2))) \to {\rm HC}^\sigma_3(\A({\rm SU}_q(2)))
\]
is surjective.
Therefore, the
Hochschild cohomology class of the Chern character is nontrivial in ${\rm HH}^3_\sigma(\A({\rm SU}_q(2)))$.
Similar comments apply to the Chern character in the next section.

\subsection[The modular Fredholm module of ${\rm SU}_q(2)$ from its spectral triple]{The modular Fredholm module of $\boldsymbol{{\rm SU}_q(2)}$ from its spectral triple}
\label{sec:the-modular-fred}

The Fredholm module for $\A({\rm SU}_q(2))$ presented above
gives (up to sign) the same (ordinary) $K$-homology class
as the Fredholm module arising from the spectral triple over $\A({\rm SU}_q(2))$
discovered in \cite{DLSSV}. It is therefore not surprising that  the modular
index pairings  and  the twisted cyclic three-cocycles obtained from these two
examples are both
nontrivial.

Let us brief\/ly recall the construction of the equivariant spectral triple
over $\A({\rm SU}_q(2))$ due to~\cite{DLSSV}. We will use the notation from that work, and refer
there for more details.

When $j \in \frac{1}{2} \N$ is a half integer, put $j^\pm=j\pm 1/2$ when $j \in \frac{1}{2} \N$.  For each $j \in \frac{1}{2} \N$, consider f\/inite-dimensional Hilbert spaces
\begin{gather*}
\Wua_j  = \mbox{span} \big\{ \ketua{j \mu n} \mid \mu = -j, -j + 1, \ldots, j, \text{ and } n = - j^+, \ldots, j^+ \big\},\\
\Wda_j  = \mbox{span} \big\{ \ketda{j \mu n} \mid \mu = -j, -j + 1, \ldots, j, \text{ and } n = - j^-, \ldots, j^- \big\},
\end{gather*}
where the elements of the respective sets form orthonormal bases.  The spectral triple is realized on the completion $\H$ of the pre-Hilbert space $\oplus_{2j = 0}^\infty \Wua_j \oplus \Wda_j$.  The action of $\A({\rm SU}_q(2))$ on $\H$ is given (\cite{DLSSV}, Proposition~4.4)
as follows.  First, the action of $a$ is given by
\begin{gather*}
\pi'(a) \ketua{j \mu n}  = \alpha^+_{j \mu n \uparrow \uparrow} \ketua{j^+ \mu^+ n^+} + \alpha^+_{j \mu n \downarrow \uparrow} \ketda{j^+ \mu^+ n^+} + \alpha^-_{j \mu n \uparrow \uparrow} \ketua{j^- \mu^+ n^+},\\
\pi'(a) \ketda{j \mu n}  = \alpha^+_{j \mu n \downarrow \downarrow} \ketda{j^+ \mu^+ n^+} + \alpha^-_{j \mu n \downarrow \downarrow} \ketda{j^- \mu^+ n^+} + \alpha^-_{j \mu n \uparrow \downarrow} \ketua{j^- \mu^+ n^+},
\end{gather*}
where the coef\/f\/icients $\alpha^{\pm}_{j \mu n}$ are given by (writing $[k]=(q^{-k}-q^k)(q^{-1}-q)^{-1}$)
\begin{gather*}
\left(\begin{matrix}
\alpha^+_{j \mu n \uparrow\uparrow} & \alpha^+_{j \mu n \uparrow\downarrow}  \\
\alpha^+_{j \mu n \downarrow\uparrow} & \alpha^+_{j \mu n \downarrow\downarrow}
\end{matrix}\right) =
 q^{\big(\mu + n - {\frac{1}{2}}\big)/2} \qn{j + \mu + 1}^{\frac{1}{2}} \left (\begin{matrix}
 q^{-j-{\frac{1}{2}}}\frac{\qn{j+n+\frac{3}{2}}^{\frac{1}{2}}}{\qn{2 j + 2}} & 0\\
 q^{\frac{1}{2}} \frac{\qn{j-n+{\frac{1}{2}}}^{\frac{1}{2}}}{\qn{2j+1}\qn{2j+2}} & q^{-j} \frac{\qn{j+n+{\frac{1}{2}}}^{\frac{1}{2}}}{\qn{2j+1}}
\end{matrix} \right )
\end{gather*}
and
\begin{gather*}
\left(\begin{matrix}
\alpha^-_{j \mu n \uparrow\uparrow} & \alpha^-_{j \mu n \uparrow\downarrow}  \\
\alpha^-_{j \mu n \downarrow\uparrow} & \alpha^-_{j \mu n \downarrow\downarrow}
\end{matrix}\right) =
q^{\big(\mu + n - {\frac{1}{2}}\big)/2} \qn{j - \mu}^{\frac{1}{2}} \left (\begin{matrix}
 q^{j+1}\frac{\qn{j-n+{\frac{1}{2}}}^{\frac{1}{2}}}{\qn{2j+1}} & - q^{\frac{1}{2}} \frac{\qn{j+n+{\frac{1}{2}}}^{\frac{1}{2}}}{\qn{2j}\qn{2j+1}}\\
 0 & q^{j+{\frac{1}{2}}}\frac{\qn{j-n-{\frac{1}{2}}}^{\frac{1}{2}}}{\qn{2j}}
\end{matrix} \right).
\end{gather*}

Similarly, the action of $b$ can be expressed as
\begin{gather*}
\pi'(b) \ketua{j \mu n}  = \beta^+_{j \mu n \uparrow \uparrow} \ketua{j^+ \mu^+ n^-} + \beta^+_{j \mu n \downarrow \uparrow} \ketda{j^+ \mu^+ n^-} + \beta^-_{j \mu n \uparrow \uparrow} \ketua{j^- \mu^+ n^-},\\
\pi'(b) \ketda{j \mu n}  = \beta^+_{j \mu n \downarrow \downarrow} \ketda{j^+ \mu^+ n^-} + \beta^-_{j \mu n \downarrow \downarrow} \ketda{j^- \mu^+ n^-} + \beta^-_{j \mu n \uparrow \downarrow} \ketua{j^- \mu^+ n^-},
\end{gather*}
where the coef\/f\/icients are given by
\begin{gather*}
\left(\begin{matrix}
\beta^+_{j \mu n \uparrow\uparrow} & \beta^+_{j \mu n \uparrow\downarrow}  \\
\beta^+_{j \mu n \downarrow\uparrow} & \beta^+_{j \mu n \downarrow\downarrow}
\end{matrix}\right) =
 q^{\big(\mu + n - {\frac{1}{2}}\big)/2} \qn{j + \mu + 1}^{\frac{1}{2}} \left (\begin{matrix}
 \frac{\qn{j-n+\frac{3}{2}}^{\frac{1}{2}}}{\qn{2 j + 2}} & 0\\
 -q^{-j-1} \frac{\qn{j+n+{\frac{1}{2}}}^{\frac{1}{2}}}{\qn{2j+1}\qn{2j+2}} & q^{-{\frac{1}{2}}} \frac{\qn{j-n+{\frac{1}{2}}}^{\frac{1}{2}}}{\qn{2j+1}}
\end{matrix} \right )
\end{gather*}
and
\begin{gather*}
\left(\begin{matrix}
\beta^-_{j \mu n \uparrow\uparrow} & \beta^-_{j \mu n \uparrow\downarrow}  \\
\beta^-_{j \mu n \downarrow\uparrow} & \beta^-_{j \mu n \downarrow\downarrow}
\end{matrix}\right) =
q^{\big(\mu + n - {\frac{1}{2}}\big)/2} \qn{j - \mu}^{\frac{1}{2}} \left (\begin{matrix}
 -q^{-{\frac{1}{2}}}\frac{\qn{j+n+{\frac{1}{2}}}^{\frac{1}{2}}}{\qn{2j+1}} & - q^j \frac{\qn{j-n+{\frac{1}{2}}}^{\frac{1}{2}}}{\qn{2j}\qn{2j+1}}\\
 0 & -\frac{\qn{j+n-{\frac{1}{2}}}^{\frac{1}{2}}}{\qn{2j}}
\end{matrix} \right).
\end{gather*}

The Dirac operator $D$ acts as the scalar $j$ on $\Wua_{j}$ and as
$-j$ on $\Wda_{j}$.  The phase $F$ of $D$ is therefore given by the factor~$1$ on~$\Wua_{j}$ and by~$-1$ on~$\Wda_{j}$.

In this basis, the modular element $K$ is represented by
\[
K\ket{j \mu n \uparrow \downarrow} = q^{-2 (\mu + n)} \ket{j \mu n\uparrow \downarrow}.
\]
We take the von Neumann algebra~$\B(\H)$, and the weight $\Phi(T):=\tr\big(K^{1/2}TK^{1/2}\big)$ for $0\leq T\in \B(\H)$.

\begin{proposition}
The triple $(\A({\rm SU}_q(2)), H, F)$ is an odd $3$-summable modular Fredholm module
with respect to $(\B(\H),\Phi)$.
\end{proposition}

\begin{proof}
Since $x\mapsto[F, x]$ is a derivation, we only need to verify the summability condition for the generators
$x = a, b$. Let $P^\uparrow$ (resp.\ $P^\downarrow$) denote the projection onto $\oplus_j \Wua_j$ (resp.\ $\oplus_j \Wda_j$).  Then the commutator $[F, x]$ can be expressed as $\Pua x \Pda - \Pda x \Pua$.  Thus, for example,
 \begin{gather}
[F, a] \ketua{j \mu n}  \mapsto q^{\big(\mu + n - {\frac{1}{2}}\big)/2} \qn{j + \mu + 1}^{\frac{1}{2}} q^{\frac{1}{2}} \frac{\qn{j-n+{\frac{1}{2}}}^{\frac{1}{2}}}{\qn{2j+1}\qn{2j+2}} \ketda{j^+ \mu^+ n^+},
\nonumber\\
[F, a] \ketda{j \mu n}  \mapsto q^{\big(\mu + n - {\frac{1}{2}}\big)/2} \qn{j - \mu }^{\frac{1}{2}} q^{\frac{1}{2}} \frac{\qn{j+n+{\frac{1}{2}}}^{\frac{1}{2}}}{\qn{2j}\qn{2j+1}} \ketua{j^- \mu^+ n^+}.
\label{eq:a2}
 \end{gather}
Therefore, we need to establish that the coef\/f\/icients in the above expressions are summable with respect to
the modular weight. The asymptotics of $\qn{k}$ is the same as that of $q^{-k}$ as $k$ tends to inf\/inity.
Hence the asymptotics of the f\/irst component of $[F,a] K^{\frac{1}{3}}$ is bounded from above by
\[
\frac{q^{-\big(j+ \frac{2}{3}\mu +\frac{1}{3}n\big)}}{q^{-4j}}.
\]
Similarly, from ~\eqref{eq:a2}, the second component of $[F,a]K^{\frac{1}{3}}$ is
asymptotically bounded from above by
\[
\frac{q^{-\big(j+\frac{1}{3}\mu - \frac{2}{3}n\big)}}{q^{-4j}},
\]
and one can see that it is a trace class operator. Analogously for $x = b$, using the expression of
the matrices $\beta^\pm_{j \mu n}$, the `matrix coef\/f\/icients' of $[F, b] K^{\frac{1}{3}}$ are
asymptotically bounded from above by
\[
\max\left(\frac{q^{-\big(2j+\frac{2}{3}\mu + \frac{2}{3}n\big)}}{q^{-4j}},
\frac{q^{-\big(j+\frac{1}{3}\mu - \frac{1}{3}n\big)}}{q^{-4j}}\right),
\]
and \looseness=1 similar analysis shows that $[F, b] K^{\frac{1}{3}} \in L^1(H) \subset L^3(H)$.
This proves the assertion. Observe that the Fredholm module
$(\A({\rm SU}_q(2)), H, F)$ is not $2$-summable, as could be
easily seen by computing the asymptotics of $[F, b] K^{\frac{1}{2}}$, which shows that this is only bounded
but not compact.
\end{proof}

Since the product of at least three commutators with $F$ is in the domain of the modular
weight $\Phi(\cdot) = \tr \big(K^{1/2} \cdot K^{1/2}\big)$ we can def\/ine the twisted Chern character of
the modular Fredholm module as before:
\[
\Ch_{3}(x_0,x_1,x_2,x_3) = \lambda_3\frac{1}{2} \Tr(F [F, x_0] [F, x_1] [F, x_2] [F, x_3] K),\qquad x_i\in \A({\rm SU}_q(2)).
\]
We now compute the pairing of $\Ch_{3}$ with the equivariant odd $K$-group. Taking $V$ as in
(\ref{eqK1}) and with the same extension of $K$ to $\H \otimes \C^2$, we obtain the following:
\begin{proposition}
\label{propqV}
The modular index of $V$ relative to the above Fredholm module is equal to $1$.
\end{proposition}

\begin{proof}
First, observe that $V$ can be written as $ V = S U$, where $S$ and $U$ are given by
\begin{gather*}
S  =
 \left(\begin{matrix}
 0 & 1 \\ 1 & 0
\end{matrix}\right),\qquad
U  =
  \left(\begin{matrix}
  a & b \\ -q b^* & a^*
\end{matrix}\right).
\end{gather*}

Recall that by \cite{SDLSV}, the operator $PUP$, where
$P = \frac{1}{2}(1+F) \otimes \id$, has a trivial cokernel,
whereas its kernel is one-dimensional and spanned by
\[
 \xi_0 = \left ( \begin{matrix}
 \ketua{0, 0, -{\frac{1}{2}}} \\
 -q^{-1} \ketua{0, 0, {\frac{1}{2}}}
 \end{matrix} \right ).
\]
Since the matrix $S$ commutes with the projection $P$, $\xi_0$ also spans
the kernel of $PVP$. It is then easy to check that the eigenvalue of the
modular operator $\tilde{K}$ acting on $\xi_0$ is $1$.
\end{proof}

\begin{remark}
The above computations show that, the two modular Fredholm modules of Sections~\ref{sec:the-basic-fred} and~\ref{sec:the-modular-fred} are related by multiplying a nontrivial $1$-dimensional character of $\mathrm{U}(1)$ if one considers the associated $\mathrm{U}(1)$-equivariant $K$-homology classes.
\end{remark}

\section{Conclusions}

The signif\/icance of the results in \cite{CC} is that we can represent those cyclic
cocycles arising from traces on $\J^n$ as Chern characters of $n$-summable semif\/inite
Fredholm modules. Theorem~\ref{thm:main-twist} shows that we can represent those
twisted cyclic cocycles arising from KMS weights on $\J^n$ as Chern characters of
modular Fredholm modules.

This Fredholm module approach to twisted traces works well, and in the examples
avoids the dimension drop
phenomena which plague $q$-deformations. An unbounded approach, in the spirit of spectral triples,
is still a work in progress, but see \cite{KS, KRS,RS}.

\subsection*{Acknowledgements}

AR was supported by the Australian Research Council, and thanks Jens Kaad for numerous discussions on related topics. MY was
supported in part by the ERC Advanced Grant 227458
OACFT ``Operator Algebras and Conformal Field Theory''. AS acknowledges support of MNII grant 189/6.PRUE/2007/7
and thanks for the warm hospitality at Mathematical Sciences Institute, Australian
National University, Canberra.

\pdfbookmark[1]{References}{ref}
\LastPageEnding

\end{document}